\documentclass[smallextended]{svjour3}
\smartqed  
\DeclareMathAlphabet{\mathpzc}{OT1}{pzc}{m}{it}
\usepackage{mathrsfs}
\usepackage{graphicx}
\usepackage{amsfonts}
\usepackage{epsfig}
\usepackage{amssymb}
\usepackage{amsmath}

\usepackage{mathtools}
\usepackage{relsize}
\usepackage{xcolor}
\usepackage{algorithmic}
\usepackage[ruled,vlined]{algorithm2e}

\numberwithin{stheorem}{section}

\numberwithin{sremark}{section}

\def\tfrac#1#2{{\frac{\lower.6ex
\hbox{$\scriptstyle#1$}}
{\raise.7ex
\hbox{$\scriptstyle#2$}}}}

\def\NN{{\mathbb N}}
\def\RR{{\mathbb R}}

\def\Frac#1#2{\displaystyle{\frac{#1}{#2}}}

\numberwithin{equation}{section}

\begin{document}

\title{ A numerical algorithm for computing the zeros of parabolic cylinder functions in the complex plane
}

\titlerunning{  A numerical algorithm for computing the zeros of parabolic cylinder functions   }    

\author{T. M. Dunster  \and A. Gil \and \\
 D. Ruiz-Antol\'{\i}n \and J. Segura
}

\authorrunning{T.M. Dunster, A. Gil, D. Ruiz-Antol\'{\i}n, J. Segura}

\institute{T. M. Dunster \at
              Department of Mathematics and Statistics\\
              San Diego State University, 5500 Campanile Drive San Diego, CA, USA\\
              \email{mdunster@sdsu.edu}   
           \and
           A. Gil \at
              Departamento de Matem\'atica Aplicada y CC. de la Computaci\'on\\
 Universidad de Cantabria, 39005-Santander, Spain\\
              \email{amparo.gil@unican.es} 
\and
      D. Ruiz-Antol\'{\i}n \at
              Departamento de Matem\'atica Aplicada y CC. de la Computaci\'on\\
Universidad de Cantabria, 39005-Santander, Spain\\
              \email{diego.ruizantolin@unican.es} 
\and
           J. Segura \at
              Departamento de Matem\'aticas, Estad\'{\i}stica y Computaci\'on\\
            Universidad de Cantabria, 39005-Santander, Spain\\
              \email{segurajj@unican.es}          
}

\date{Submitted: Nov 27, 2024 }

\maketitle

\begin{abstract}
A numerical algorithm (implemented in Matlab) for computing the zeros of the parabolic cylinder function $U(a,z)$ in domains of the complex plane is presented.
The algorithm uses accurate approximations to the first zero plus a highly efficient method based on a fourth-order fixed point method  
with the parabolic cylinder functions computed
by Taylor series and carefully selected steps, to compute the 
rest of the zeros. 
For $|a|$ small, the asymptotic approximations are complemented with a few fixed point iterations requiring  the evaluation of $U(a,z)$ and $U'(a,z)$
in the region where the complex zeros are located.
Liouville-Green expansions are derived to enhance the performance of a computational scheme to evaluate $U(a,z)$ and $U'(a,z)$ in that region. 
Several tests show the accuracy and efficiency of the numerical algorithm. 

\keywords{Parabolic cylinder functions \and Complex zeros \and Fixed point methods \and Matlab software.      }
 \subclass{65H05 \and   33F05\and 34C10 \and 30E10      
       }
\end{abstract}

\section{Introduction} 

The search for complex zeros of functions in a domain of the complex plane is essential in various fields, particularly in problems related to the propagation and radiation of electromagnetic waves. In such contexts, the complex zeros of specific functions—such as wave propagation equations \cite{Kwon:2004:MWM}, Green’s functions, or impedance functions—play a critical role in determining key physical phenomena. These include resonance, scattering, and the stability of wave solutions. In general, the complexity of these problems often requires numerical techniques and sophisticated algorithms for accurately identifying the zeros in a given domain of the complex plane; see for example 
\cite{Michal:2018:CBT,Kowal:2018:CRP,Chen:2017:RLO}. 

Parabolic cylinder functions arise naturally as solutions to the wave equation when expressed in a parabolic coordinate system. In this paper, we present a numerical algorithm for finding the complex zeros of the parabolic cylinder function  $U(a,z)$ in a domain of the complex plane.
The function $U(a,z)$ is a solution of 
the homogeneous equation

\begin{equation}  \label{eq:01}
\frac{d^{2}y}{dz^{2}}-\left(\frac{1}{4}z^{2}+a \right)y=0.
\end{equation}

A Poincar\'e-type expansion for this function is given by  \cite[Eq. 12.9.1]{NIST:DLMF} 

\begin{equation}  \label{eq:02}
U\left(a,z\right)\sim e^{-\frac{1}{4}z^{2}}z^{-a-\frac{1}{2}}\sum_{s=0}^{%
\infty}(-1)^{s}\frac{{\left(\frac{1}{2}+a\right)_{2s}}}{s!(2z^{2})^{s}},\,\,|\arg(z)|<3\pi/4,
\end{equation}
from which, it is clear its recessive behavior at infinity in the sector $|\arg(z)| \leq \pi/4$.
 
The  algorithm uses as starting values the asymptotic approximations to the zeros of $U(a,z)$ given in \cite{Dunster:2024:UAZ}. 
These asymptotic approximations are highly accurate for moderate to large values of $a$. 
The approximations are expressed in terms of the zeros of Airy functions or combinations of these functions. For small values of 
$a$, the asymptotic approximations are refined using a fourth-order fixed-point method. The remaining zeros of 
$U(a,z)$ are obtained through a highly efficient scheme involving carefully selected steps and additional fixed-point iterations
with the parabolic cylinder functions computed
by Taylor series. 
Numerical tests demonstrate the accuracy and efficiency of the numerical scheme.
The algorithm represents a first practical implementation that illustrates how combining asymptotic 
and iterative methods is a highly efficient strategy for determining zeros of functions that are solutions to second-order ODEs.

\section{Algorithm for computing the complex zeros of $U(a,z)$}

We follow the results given in  \cite{Segura:2013:CCZ}.
In this reference, the complex zeros of solutions of ODEs

\begin{equation}
 y''(z)+A(z)y(z) = 0
\end{equation}
 with $A(z)$ a complex meromorphic function, are considered. 
It is shown that the zeros
lie over certain curves which
follow very closely 
the approximate anti-Stokes lines (ASLs). 
 Information on the approximate Stokes lines (SLs) is also important
in a general strategy for finding the complex zeros of $y(z)$.
This qualitative analysis of ASLs and SLs is then combined with the use of a fixed point iteration
 $T(z)=z+\frac{1}{\sqrt{A(z)}} \arctan  \left(  \sqrt{A(z)}  \frac{y(z)}{y'(z)} \right)$ and carefully selected  displacements 
$H^{\pm} (z)=   z \pm \frac{\pi}{\sqrt{A(z)}} $.  

 The strategy for finding the complex zeros of $y(z)$  can be summarized as follows:

\begin{enumerate}
\item Divide the complex plane in disjoint domains separated by the principal ASLs
and SLs and compute separately in each domain.  Schwarz symmetry can be used to reduce the problem.

\item In each domain, start away from the principal SLs, close to a principal ASL and/or
singularity (if any). Iterate with $T(z)$ until a first zero is found or use an asymptotic approximation
(if available) for that zero.  

\item Then, use the basic algorithm described in \cite{Segura:2013:CCZ} for computing consecutive zeros, 
choosing the displacements $H^{\pm}(z)$ in the
direction of approaching the principal SLs and/or singularity.
\end{enumerate}

For the parabolic cylinder function $U(a,z)$, we consider real orders $a$ (excluding the case $a = -k+\frac12,\,k \in \NN$, which corresponds to
the case of Hermite polynomials). For $a \in \RR, \,a \neq -k+\frac12,\,k \in \NN$, there are an infinite number of complex zeros of the function $U(a,z)$
 tending to the ray $\arg z = 3\pi/4$ and a conjugate string. 
To compute the complex zeros of $U(a,z)$ with $a<0$, the following displacements  $H^{+}(z)$ and iterating function $T(z)$ 
are used:

\begin{equation}
H^+(z)=z+\Frac{\pi}{\sqrt{-z^2/4-a}}\,,
\label{displa}
\end{equation}

\begin{equation}
T(z)=z-\Frac{1}{\sqrt{-z^2/4-a}}\arctan\left(\sqrt{-z^2/4-a} \,Q(a,z)\right)\,,
\label{iter}
\end{equation}
where 

\begin{equation}
Q(a,z)=\Frac{U(a,z)}{U'(a,z)}\,.
\label{quot}
\end{equation}

An algorithm to compute the complex zeros of $U(a,z)$ for  $a<0$,  in the domain of the complex plane $\Im z \in [0,\,L]$, $\Re z<0$ is described in Algorithm \ref{algo}.
To compute the zeros for $a>0$,  a similar algorithm can be used but considering the function $U(a,iz)$.
 Also, the factor $\displaystyle\sqrt{-z^2/4-a}$
appearing in (\ref{displa}) and (\ref{iter})
should be changed to $ \displaystyle\sqrt{-z^2/4+a}$. For $a>0$, the complex zeros are obtained in the domain of the complex plane $\Re z \in [-L,\,0]$, $\Im z>0$.

In Algorithm \ref{algo}, we consider $z=-L+iL$ as a starting  point. Then, the zero $z_m$ of $U(a,z)$ closest to $z$ is calculated using the asymptotic approximations given in 
   \cite{Dunster:2024:UAZ}. The value of $m$ is estimated using the first term in the approximation  \cite[Eq. 12.11.1]{NIST:DLMF} 

\begin{equation}
z \approx e^{3/4 \pi i} \sqrt{2 \tau_m},
\end{equation}
\hspace{-0.1cm}
where 
\[\tau_{m}=\left(2m+\tfrac{1}{2}-|a|\right)\pi+i\ln\left(\pi^{-\frac{1}{2}}2^{-|a|-\frac{1}{2}}\Gamma\left(\tfrac{1}{2}+|a|\right)\right).\]

For $a<0$,  the asymptotic approximations $z_m$  to the zeros of $U(a,z)$ are given in terms of the zeros of the following combination of Airy functions

\begin{equation}
\label{eq:air}
\mathcal{A}i(a,z)
=2e^{-\pi i/6} \cos\left(a \pi\right)\mathrm{Ai}_{1}(z)
+ie^{a\pi i}\mathrm{Ai}(z),
\end{equation}
where $\mathrm{Ai}_{1}(z)=\mathrm{Ai}(ze^{-2\pi i/3})$.

For $a>0$, the asymptotic approximations $z_m$ are given in terms of the negative zeros, $a_m$, of the Airy functions $\mathrm{Ai}(x)$.
In our algorithm, for $m$ moderate, we use precomputed values of these zeros $a_m$; for $m$ large, we use the expansions \cite[Eq. 9.9.6]{NIST:DLMF} 
\[a_{m}=-T\left(\tfrac{3}{8}\pi(4m-1)\right),\] where

$$
\begin{array}{lcl}
T(t) &\sim & t^{2/3}\left(1+\Frac{5}{48}t^{-2}-\Frac{5}{36}t^{-4}+\Frac{77125}{829%
44}t^{-6}-\Frac{1080\;56875}{69\;67296}t^{-8}\right.\\
 &&+\left.\Frac{16\;23755\;96875}{3344\;30%
208}t^{-10}-\cdots\right).
 \end{array}
$$

For $|a|$ small or $m$ small, the approximations to the zeros obtained using asymptotic expansions are refined with a few iterations
 of \eqref{iter}. The computation of $U(a,z)$ and $U'(a,z)$
needed in  \eqref{quot} is discussed in Section \ref{compU}.  

Using the previous approximations, the first zero is calculated $z_c^{(0)}=z_m$ in Algorithm \ref{algo}.
For computing the second zero $z_c^{(1)}$, the  step $h=\Frac{\pi}{\sqrt{-z^2/4-a}}$ in \eqref{displa} is taken. Then, we consider $z=z_c^{(0)}+h$ and evaluate $z=T(z)$
using Taylor series centered at $z_c^{(0)}$ to compute \eqref{quot}. For computing the Taylor series, we will use the results given in Section  \ref{taylor}; in the first iteration, 
we use that $U(a, z_c^{(0)})=0$, and for the derivative, since we are interested in the zeros, it is possible to take a normalized value ($U'(a, z_c^{(0)})=1$,
for example).  The rest of fixed point iterations in the inner while loop in Algorithm \ref{algo} are also computed using Taylor series.
The algorithm stops when a zero with an imaginary part smaller than $\delta$ ($\delta$ being small and positive) is computed.

\begin{algorithm}[h!]
\SetAlgoLined
\KwData{$a$, real negative parameter; $L$, length of the interval.}
\KwResult{complex zeros $z_c^{(j)},\,j=0,1,2,...$ in $\Im z \in [0,\,L]$, $\Re z<0$. }
\begin{enumerate} 

\item Set $z=-L+iL$; $\epsilon=10^{-14}$; $\delta=10^{-4}$.

\item Calculate $z_c^{(0)}$ (closest zero to $z$) using an asymptotic expansion in terms
of the zeros of \eqref{eq:air}; 

\item If $a$ is small, refine the value iterating $T(z)$. 

\item Set $U(a,z_c^{(0)})=0$, $U'(a,z_c^{(0)})=1$ (function values for the Taylor series);

\item $i=0$\;

\item

\hspace*{-0.2cm} \While{ $\Im z_c^{(i)}>\delta$}{
 $z=H^+(z_c^{(i)}); \Delta=1+\epsilon$\;
 \While{$\Delta > \epsilon$ }{
   $y=z$; $z=T(z)$;

  $\Delta = |z-y|/|y|$\;
  }
 $i=i+1$;

$z^{(i)}_c=z$.

Set $U(a,z_c^{(i)})=0$, $U'(a,z_c^{(i)})=1$ (function values for the Taylor series);

}
\item Check that all zeros $z_c^{(j)},\,j=0,1,2,...$ satisfy $\Im z_c^{(j)}  \in [0,\,L]$.
\end{enumerate}
\label{algo}
 \caption{Computation of the complex zeros of $U(a,z),\,a<0$ in the domain
 $\Im z \in [0,\,L]$, $\Re z<0$.}
\end{algorithm}

\section{Local Taylor series}
\label{taylor}

Assuming that $y(z_0)$ and $y'(z_0)$ are available, it is possible to compute
the functions  $y(z_1)=y(z_0+h)$ and $y'(z_1)=y'(z_0+h)$ using Taylor series for $y(z)$ and $y'(z)$ around
$z_0$. That is, we
compute

\begin{equation}
\begin{array}{l}
\label{Taylorsums}
y (z_{i+1})=\displaystyle\sum_{k=0}^{N}y^{(k)}(z_i)\Frac{h^k}{k!}+  {\cal O} (h^{N+1}),\\
y' (z_{i+1})=\displaystyle\sum_{k=0}^{N}y^{(k+1)}(z_i)\Frac{h^k}{k!}+  {\cal O} (h^{N+1}).
\end{array}
\end{equation}


The successive derivatives $y^{(k)}$ can be computed by differentiation of the differential 
equation. From \eqref{eq:01} we have, differentiating $k$ times, $k\ge 2$:
\begin{equation}
\label{recuder}
y^{(k+2)}-\left(\tfrac14z^2+a\right)y^{(k)}-\tfrac12zky^{(k-1)}-\tfrac14k(k-1)y^{(k-2)}=0,
\end{equation}
which allows computing derivatives at $z=z_i$  when  $y^{(j)}(z_i)$, $j=0,1,2,3$ are known.
$y^{(0)}(z_i)$ and $y^{(1)}(z_i)$ are known from the previous step, and
\begin{equation}
\begin{array}{l}
y^{(2)}(z_i)=\left(\tfrac14 z_i^2+a\right)y^{(0)}(z_i),\\
y^{(3)}(z_i)=\left(\tfrac14 z_i^2+a\right)y^{(1)}(z_i)+\tfrac12  z_i y^{(0)}(z_i).
\end{array}
\end{equation}

 For computing the successive derivatives, it is necessary that the recursion process 
 for Eq. (\ref{recuder}) is well conditioned. 
Using the Perron-Kreuser theorem \cite{Perron:1914:HRR},
 it is easy to see that there are no exponentially dominant nor recessive solutions 
of the linear recurrence relation  \eqref{recuder}. Forward computation, therefore, is not seriously ill conditioned as $k$
 becomes large.

\section{Computation of $U(a,z)$ and its derivative in the region where the complex zeros lie}
\label{compU}

In \cite{Dunster:2024:PCF}, methods to compute $U(a,z)$ were given. A computational scheme based on these methods (Airy-type expansions, integral representations, Maclaurin series, and Poincar\'e expansions) can be designed to evaluate $U(a,z)$.  For the derivative $U'(a,z)$, the relation \cite[Eq. 12.8.3]{NIST:DLMF}

\[U'\left(a,z\right)-\tfrac{1}{2}zU\left(a,z\right)+U\left(a-1,z\right)=0,\]
could be used. 

In the region where the complex zeros of the function are located (the second and third quadrants of the complex plane), the scheme can be further enhanced by employing Liouville-Green approximations, which we discuss next.

\subsection{Liouville-Green expansions}

As in \cite[Eq. (2.2)]{Dunster:2021:UAP} we define a Liouville-Green variable $\bar{\xi}$ given  by
\begin{equation}
\label{26a}
\bar{\xi}= \int_{0}^{\hat{z}}\left(t^2+1\right)^{1/2} dt
=\frac{1}{2}\hat{z}\left(  \hat{z}^{2}+1 \right)^{1/2}   +\frac{1}{2}\ln\left( \hat{z}+\left(\hat{z}^{2}+1  \right)^{1/2}   \right),
\end{equation}
where here and throughout bars do not denote complex conjugate, unless otherwise noted. The branch is chosen so that $\bar{\xi}$ is real when $\hat{z}$ is, both being of the same sign, and by continuity elsewhere in the plane with cuts along $\hat{z}=\pm iy$, $1\leq y < \infty$. We are only interested for $\hat{z}$ in the second quadrant, since the complex zeros of $U(a,z)$ lie there (there is also a conjugated set of zeros on the third quadrant).

Thus we find that as $\hat{z} \rightarrow \infty$ with $\Re(\hat{z}) < 0$
\begin{equation}
\label{26b}
\bar{\xi}=
-\tfrac{1}{2}\hat{z}^{2} - \tfrac{1}{2}\ln(-2\hat{z})-\tfrac{1}{4}
+\mathcal{O}(\hat{z}^{-2}).
\end{equation}
In particular $\bar{\xi} \rightarrow -\infty$ as $\hat{z} \rightarrow -\infty$ and $\Re(\bar{\xi}) \rightarrow +\infty$ as $\hat{z} \rightarrow i\infty$ to the left of the cut.

Next from \cite[Eq. (2.8)]{Dunster:2021:UAP} let
\begin{equation}
\label{26c}
\bar{\beta}=\frac{\hat{z}}{\sqrt{\hat{z}^{2}+1}},
\end{equation}
with $\bar{\beta}<0$ when $\hat{z}<0$ and is continuous in the same cut plane as for $\bar{\xi}$. Thus $\bar{\beta} \rightarrow -1$ as $\hat{z} \rightarrow \infty$ in the left half plane $\Re(\hat{z}) < 0$.

Then from (2.14), (2.16)-(2.18), (2.20) and (3.10) of that paper as $u \rightarrow \infty$ we find the expansion
\begin{multline}
\label{26g}
U\left(\tfrac{1}{2}u,-\sqrt {2u} \, \hat{z}\right)
\sim \left( \frac{2e}{u} \right)^{u/4}
\frac{1}{\left\{2u\left(1+\hat{z}^{2}\right)
\right\}^{1/4}}
\\   \times
\exp \left\{u\bar{\xi}
+\sum\limits_{s=1}^{\infty}{(-1)^{s}
\frac{\mathrm{E}_{s}(\bar{\beta})
-\mathrm{E}_{s}(-1)}{u^{s}}}
\right\},
\end{multline}
for $\hat{z}$ in a domain that certainly contains the second quadrant, except for a closed neighborhood of the turning point $\hat{z}=i$.

The coefficients $\mathrm{E}_{s}$ are defined as

\begin{equation}
\label{07}
\mathrm{E}_{1}(\bar\beta)=\tfrac{1}{24}\bar\beta
\left(5\bar\beta^{2}-6\right),
\end{equation}
\begin{equation}
\label{08}
\mathrm{E}_{2}(\bar\beta)=
\tfrac{1}{16}\left(\bar\beta^{2}-1\right)^{2} 
\left(5\bar\beta^{2}-2\right),
\end{equation}
and for $s=2,3,4\cdots$
\begin{equation}
\label{09}
\mathrm{E}_{s+1}(\bar\beta) =
\frac{1}{2} \left(\bar\beta^{2}-1 \right)^{2}\mathrm{E}_{s}^{\prime}(\bar\beta)
+\frac{1}{2}\int_{\sigma(s)}^{\bar\beta}
\left(p^{2}-1 \right)^{2}
\sum\limits_{j=1}^{s-1}
\mathrm{E}_{j}^{\prime}(p)
\mathrm{E}_{s-j}^{\prime}(p) dp,
\end{equation}
where $\sigma(s)=1$ for $s$ odd and $\sigma(s)=0$ for $s$ even. We remark that ${\mathrm{E}}_{2s}(\bar{\beta})$ is even ${\mathrm{E}}_{2s+1}(\bar{\beta})$ is odd, and ${\mathrm{E}}_{2s}(\pm 1)=0$.

Next, from  (\ref{eq:02}) we have

\begin{equation}
\label{26h}
U\left(-\tfrac{1}{2}u,-i\sqrt {2u} \, \hat{z}\right)
\sim 
(2u)^{\frac{1}{4}u-\frac{1}{4}}
\hat{z}^{\frac{1}{2}u-\frac{1}{2}}
\exp\left\{\frac{1}{4}(1-u) \pi i
+\frac{1}{2} u \hat{z}^{2}\right\},
\end{equation}
as $\hat{z} \rightarrow i\infty$, and as such this function is recessive at $\hat{z}=i\infty$. Now again in a domain that contains the second quadrant, but this time excluding the points $\hat{z}=iy$, $0\leq y \leq 1$, we can show in a similar manner by matching recessive solutions that as $u \rightarrow \infty$
\begin{multline} 
\label{26i}
U\left(-\tfrac{1}{2}u,-i\sqrt {2u} \, \hat{z}\right)
\sim \left( \frac{u}{2e} \right)^{u/4}
\frac{e^{(u-1)\pi i/4}}
{\left\{2u\left(1+\hat{z}^{2}\right)
\right\}^{1/4}}    \\   \times
\exp \left\{-u\bar{\xi}
+\sum\limits_{s=1}^{\infty}{
\frac{\mathrm{E}_{s}(\bar{\beta})
-\mathrm{E}_{s}(-1)}{u^{s}}}
\right\}.
\end{multline}

Next from \cite[Eqs. (3.20) and (3.24)]{Dunster:2021:UAP} we have for large $u$
\begin{equation} 
\label{26j}
\frac{\sqrt {2\pi}}
{\Gamma\left(\frac{1}{2}u +\frac{1}{2} \right)}
\left(\frac{u}{2e}\right)^{u/2}
\sim 
\exp\left\{2\sum\limits_{s=0}^{\infty}
\frac{\mathrm{E}_{2s+1}(-1)}{u^{2s+1}} \right\}.
\end{equation}

Then, we have
\begin{multline} 
\label{26n}
U\left(\tfrac{1}{2}u,\sqrt {2u} \, \hat{z}\right)
\sim    \left( \frac{2e}{u} \right)^{u/4}
\frac{2e^{-(u+1)\pi i/4}}
{\left\{2u\left(1+\hat{z}^{2}\right)
\right\}^{1/4}}
\\   \times
\exp \left\{\sum\limits_{s=1}^{\infty}
\frac{\mathrm{E}_{2s}(\bar{\beta})}{u^{2s}}
+\sum\limits_{s=0}^{\infty}
\frac{\mathrm{E}_{2s+1}(-1)}{u^{2s+1}}
\right\}  \cos(\chi(u,\hat{z})),
\end{multline}
where
\begin{equation}
\label{26p}
\chi(u,\hat{z})
= iu\bar{\xi} + \frac{1}{4}(u+1)\pi
-i\sum\limits_{s=0}^{\infty}
\frac{\mathrm{E}_{2s+1}(\bar{\beta})}{u^{2s+1}}.
\end{equation}
It is evident that for large $u$ the zeros asymptotically lie on the curve $\Im\{\chi(u,\hat{z})\}=0$. Note from (\ref{26a}) that when $\hat{z}=i$
\begin{equation}
\label{26q}
i\bar{\xi}= i\int_{0}^{i}\left(t^2+1\right)^{1/2} dt
=-\frac{1}{4}\pi,
\end{equation}
and hence (\ref{26p}) can be rewritten as
\begin{equation}
\label{26r}
\chi(u,\hat{z})
= u\rho + \frac{1}{4}\pi
-i\sum\limits_{s=0}^{\infty}
\frac{\mathrm{E}_{2s+1}(\bar{\beta})}{u^{2s+1}},
\end{equation}
where
\begin{equation}
\label{26s}
\rho= i\int_{i}^{\hat{z}}\left(t^2+1\right)^{1/2} dt.
\end{equation}

The expansion (\ref{26n}) is valid in a domain which contains the second quadrant, except for points close to $\hat{z}=iy$, $0\leq y \leq 1$. In particular it is valid on the zero curve which is close to $\Re(\bar{\xi})=\mathrm{constant}$ in the second quadrant emanating from $\hat{z}=i$, except for points close to this turning point; see \cite[Fig. 1]{Dunster:2021:UAP}.

Next consider the derivative. We find that $w(\hat{z})=(\hat{z}^2+1)^{-1/2} U'(\tfrac{1}{2}u,\sqrt {2u}\, \hat{z})$ satisfies
\begin{equation}
\label{27a}
\frac{d^2w}{d\hat{z}^2} 
= \left\{u^2 \left(\hat{z}^2+1\right)
+\frac{2\hat{z}^{2}-1}
{\left(\hat{z}^{2}+1\right)^{2}}\right\}w.
\end{equation}
We follow \cite[Sect. 2]{Dunster:2021:UAP} with

\begin{equation}
\label{27b}
\Phi(\hat{z})=\frac{5\hat{z}^{2}-2}
{4\left(\hat{z}^{2}+1\right)^{3}}.
\end{equation}
The coefficients in our expansions are given by
\begin{equation} 
\label{27c}
\tilde{\mathrm{E}}_{1}(\bar{\beta})=\tfrac{1}{24}
\bar{\beta}
\left(7\bar{\beta}^{2}-6\right),
\end{equation}
and
\begin{equation} 
\label{27d}
\tilde{\mathrm{E}}_{2}(\bar{\beta})=
\tfrac{1}{16}\left(1-\bar{\beta}^{2}\right)^{2} 
\left(2-7\bar{\beta}^{2}\right),
\end{equation}
and for $s=2,3,4\cdots$
\begin{equation} 
\label{27e}
\tilde{\mathrm{E}}_{s+1}(\bar{\beta}) =
-\frac{1}{2} \left(1-\bar{\beta}^{2} \right)^{2}\tilde{\mathrm{E}}_{s}^{\prime}(\bar{\beta})
-\frac{1}{2}\int_{\sigma(s)}^{\bar{\beta}}
\left(1-p^{2} \right)^{2}
\sum\limits_{j=1}^{s-1}
\tilde{\mathrm{E}}_{j}^{\prime}(p)
\tilde{\mathrm{E}}_{s-j}^{\prime}(p) dp.
\end{equation}
Again $\sigma(s)=1$ for $s$ odd and $\sigma(s)=0$ for $s$ even, so that $\tilde{\mathrm{E}}_{2s}(\bar{\beta})$ is even, $\tilde{\mathrm{E}}_{2s+1}(\bar{\beta})$ is odd, and $\tilde{\mathrm{E}}_{2s}(\pm 1)=0$.

Then from \cite[Eq.(2.27)]{Dunster:2021:UAP} we have
\begin{multline} 
\label{27f}
U'\left(\tfrac{1}{2}u,-\sqrt {2u}\, \hat{z}\right) \sim
-\frac{1}{2}
\left(\frac{2e}{u} \right)^{u/4}
\left\{2u\left(1+\hat{z}^{2}\right)
\right\}^{1/4} \\
\times
\exp \left\{u\bar{\xi}
+\sum\limits_{s=1}^{\infty}
\frac{
\tilde{\mathrm{E}}_{s}(\bar{\beta}) 
-\tilde{\mathrm{E}}_{s}(-1)}{u^{s}}
\right\},
\end{multline}
as $u \to \infty$ for $\hat{z}$ lying in a domain that includes the second quadrant, bar $\hat{z}=i$.

Next, similarly to (\ref{26i}) we obtain for the solution of (\ref{27a}) that is recessive at $\hat{z}=i \infty$
\begin{multline} 
\label{27g}
U'\left(-\tfrac{1}{2}u,-i\sqrt {2u}\, \hat{z}\right) \sim
-\frac{1}{2}
\left( \frac{u}{2e} \right)^{u/4}
e^{(u+1)\pi i/4}
\left\{2u\left(1+\hat{z}^{2}\right)
\right\}^{1/4}
\\
\times
\exp \left\{-u\bar{\xi}
+\sum\limits_{s=1}^{n-1}(-1)^{s}
\frac{\tilde{\mathrm{E}}_{s}(\bar{\beta}) 
-\tilde{\mathrm{E}}_{s}(-1)}{u^{s}}
\right\}
\end{multline}
as $u \to \infty$. This too is valid in the second quadrant, except for points close to $\hat{z}=iy$, $0\leq y \leq 1$.

Now as $\hat{z} \to +\infty$ with $a +\frac{1}{2} \neq 0,-1,-2,\ldots$
\begin{equation}
\label{27h}
U'(a,-\hat{z})\sim 
-\sqrt {\frac{\pi}{2}}
\frac{\hat{z}^{a+\frac{1}{2}}
e^{\frac{1}{4}\hat{z}^{2}}}
{\Gamma\left(a +\frac{1}{2} \right)}.
\end{equation}
Further from (\ref{26a}) we use that as $\hat{z} \to +\infty$
\begin{equation}
\label{27i}
\bar{\xi}=
\tfrac{1}{2}\hat{z}^{2} + \tfrac{1}{2}\ln(2\hat{z})+\tfrac{1}{4}
+\mathcal{O}(\hat{z}^{-2}).
\end{equation}

On comparing (\ref{27f}) with (\ref{27h}) as $\hat{z} \to +\infty$, noting that $\beta \to 1$ in this case, and that $\tilde{\mathrm{E}}_{s}(-1)=(-1)^{s}\tilde{\mathrm{E}}_{s}(1)$, we derive similarly to (\ref{26j})

\begin{equation} 
\label{27j}
\frac{\sqrt {2\pi}}
{\Gamma\left(\frac{1}{2}u +\frac{1}{2} \right)}
\left(\frac{u}{2e}\right)^{u/2}
\sim 
\exp\left\{2\sum\limits_{s=0}^{\infty}
\frac{\tilde{\mathrm{E}}_{2s+1}(1)}{u^{2s+1}} \right\}
\quad (u \to \infty).
\end{equation}
Next, we shall use the differentiated form of \cite[Eq. 12.2.17]{NIST:DLMF} 
\begin{equation}
\label{27k}
U'(a,\hat{z}) =
ie^{-a \pi i} U'(a,-\hat{z})
+\frac{\sqrt{2\pi}\,
e^{-(\frac{1}{2}a+\frac{1}{4})\pi i}}
{\Gamma\left(a+\tfrac{1}{2}\right)}
U'(-a, -i \hat{z}).
\end{equation}
Hence from (\ref{27f}), (\ref{27g}), (\ref{27j}) and (\ref{27k}) we arrive at our desired expansion, valid as $u \to \infty$ and (at least) $\hat{z}$ lying in the second quadrant (excluding the interval $\hat{z}=iy$, $0\leq y \leq 1$)
\begin{multline} 
\label{27l}
U'\left(\tfrac{1}{2}u,\sqrt {2u} \, \hat{z}\right)
\sim    -\left( \frac{2e}{u} \right)^{u/4}
e^{-(u-1)\pi i/4}
\left\{2u\left(1+\hat{z}^{2}\right)
\right\}^{1/4}
\\   \times
\exp \left\{\sum\limits_{s=1}^{\infty}
\frac{\tilde{\mathrm{E}}_{2s}(\bar{\beta})}{u^{2s}}
+\sum\limits_{s=0}^{\infty}
\frac{\tilde{\mathrm{E}}_{2s+1}(1)}{u^{2s+1}}
\right\}  \sin(\tilde{\chi}(u,\hat{z})),
\end{multline}
where
\begin{equation}
\label{27m}
\tilde{\chi}(u,\hat{z})
= iu\bar{\xi} + \frac{1}{4}(u+1)\pi
+i\sum\limits_{s=0}^{\infty}
\frac{\tilde{\mathrm{E}}_{2s+1}
(\bar{\beta})}{u^{2s+1}}
= u\rho + \frac{1}{4}\pi
+i\sum\limits_{s=0}^{\infty}
\frac{\tilde{\mathrm{E}}_{2s+1}
(\bar{\beta})}{u^{2s+1}},
\end{equation}
with $\rho$ given by (\ref{26s}).

\begin{figure}[h!]
		\includegraphics[width=1.1\textwidth]{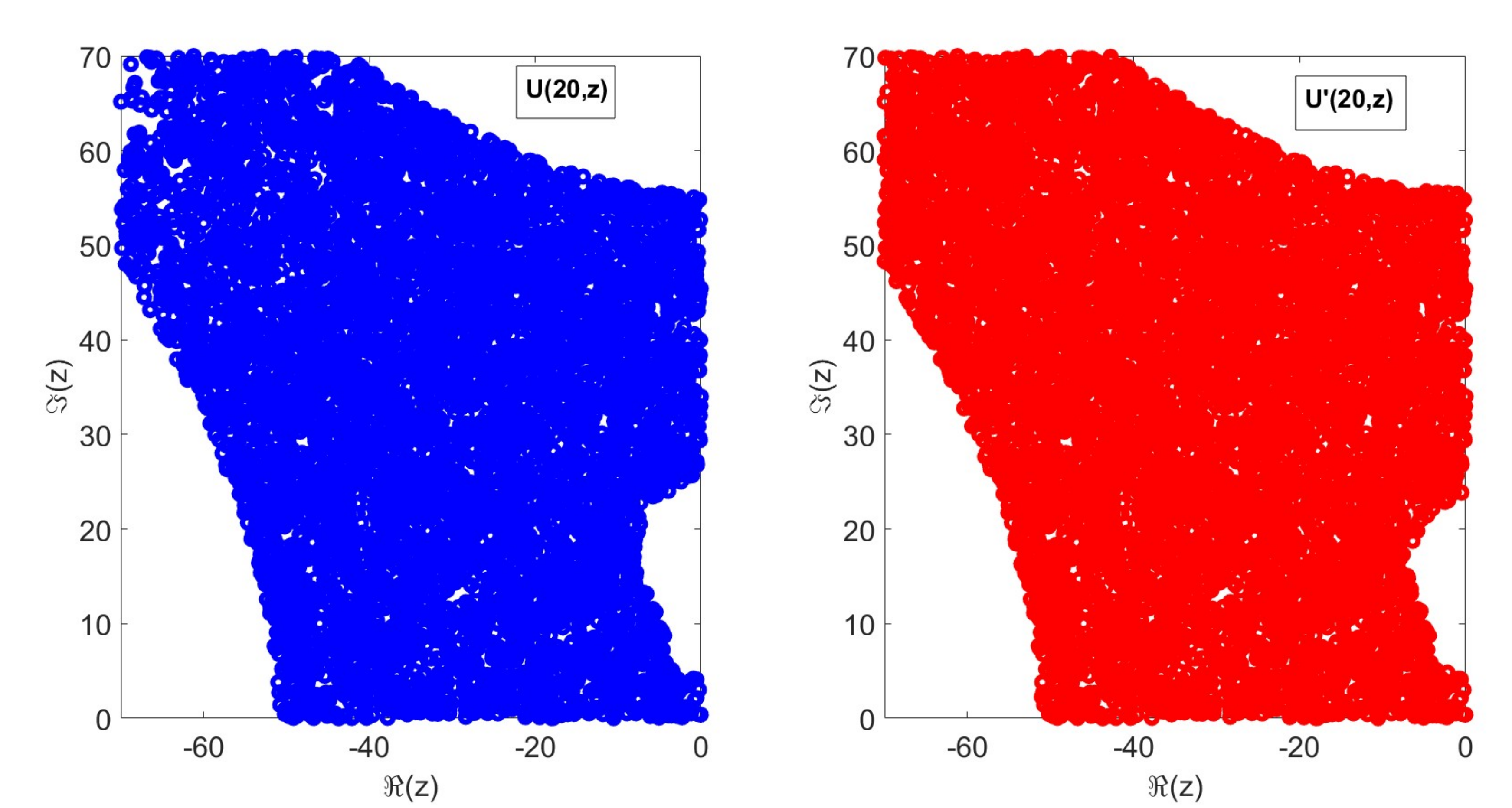}
	\caption{Test of the accuracy obtained with the Liouville-Green expansions. Left: Points where the error
when testing the recurrence relation \cite[Eq.12.8.1]{NIST:DLMF} for $U(20,z)$ is smaller than $5\times 10^{-13}$.  
 Right: Points where the error
when testing the recurrence relations    \cite[Eq.12.8.2]{NIST:DLMF} and \cite[Eq.12.8.3]{NIST:DLMF} for $U'(20,z)$ is smaller than $5\times 10^{-13}$.  }
	\label{fig:0}
\end{figure}

An example of the accuracy obtained with the use of Liouville-Green expansions for computing
$U(a,z)$ and its derivative, is given in Figure \ref{fig:0}. The Liouville-Green expansion \eqref{26n}
has been tested using the
the recurrence relation
\cite[Eq.12.8.1]{NIST:DLMF}. A large number ($10^4$) of  $z= \Re z + i\Im z$ points have been randomly generated in the domain $\Re z \in (-70,\,0)$, $\Im z \in (0,\,70)$. 
The points where the error when testing the recurrence relation for $U(20,z)$ is smaller than $5\times 10^{-13}$, are plotted in Figure \ref{fig:0} (left). 
A similar plot for the derivative $U'(a,z)$ is shown in  Figure \ref{fig:0} (right).
For testing the accuracy of the Liouville-Green expansion for the derivative \eqref{27l}, we use the relations
       \cite[Eq.12.8.2]{NIST:DLMF} and \cite[Eq.12.8.3]{NIST:DLMF}.  Our tests show that for moderate or large values of the parameter $a$, the Liouville-Green expansions allow for calculating $U(a,z)$ and $U'(a,z)$ with very high accuracy in the region where the zeros are located. In our algorithm, we use the expansions for values of $a$ greater than $18$
and $z=\Re z + i\Im z$, with $|\Re z| > 15$, $|\Im z| >15$.

\section{Numerical computation of the zeros}

The algorithms to compute the complex zeros of $U(a,z)$ for $a$ positive and negative have been implemented in Matlab.
The following Matlab function 

\vspace*{0.2cm}
\verb|[zcer]=zerosUaz(a,L)|;  
\vspace*{0.2cm}

\hspace*{-0.5cm}computes the complex zeros of $U(a,z)$  in the domain of the complex plane $\Im z \in [0,\,L]$, $\Re z<0$ (for $a<0$) and
$\Re z \in [-L,\,0]$, $\Im z>0$ (for $a>0$). \verb|zerosUaz(a,L)| makes use of
the  functions \verb|zerosUazpos(a,L)|  for $a>0$ and of the function \verb|zerosUazneg(a,L)|  for $a<0$.  
The functions \verb|zerosUaz(a,L)| can be obtained from GitHub\footnote{{https://github.com/AmparoGil/NumerZerosPCFs}}.

For testing the accuracy of the zeros  obtained with numerical algorithm $z_c^{(j)},\,j=0,1,2,...$, the relative error in the computations can be estimated using the inverse of the condition number 

\begin{equation}
\mbox{Relative error} (z_c^{(j)})\approx  \left| \frac{1}{ z_c^{(j)}  } \frac{U(a,z_c^{(j)} )}{U'(a,z_c^{(j)}) } \right|,\, j=0,1,2,...
\label{ercond}
\end{equation}

For calculating the function and its derivative, we use the results mentioned in Section 3. It is important to note that, except for the first zero (in the case where $a$ is small), the calculation of the relative error \eqref{ercond} can be considered as an independent test of the algorithm's accuracy, considering that in the fixed-point iterations for the remaining zeros, we use Taylor series.

\begin{figure}[h!]
		\includegraphics[width=1.2\textwidth]{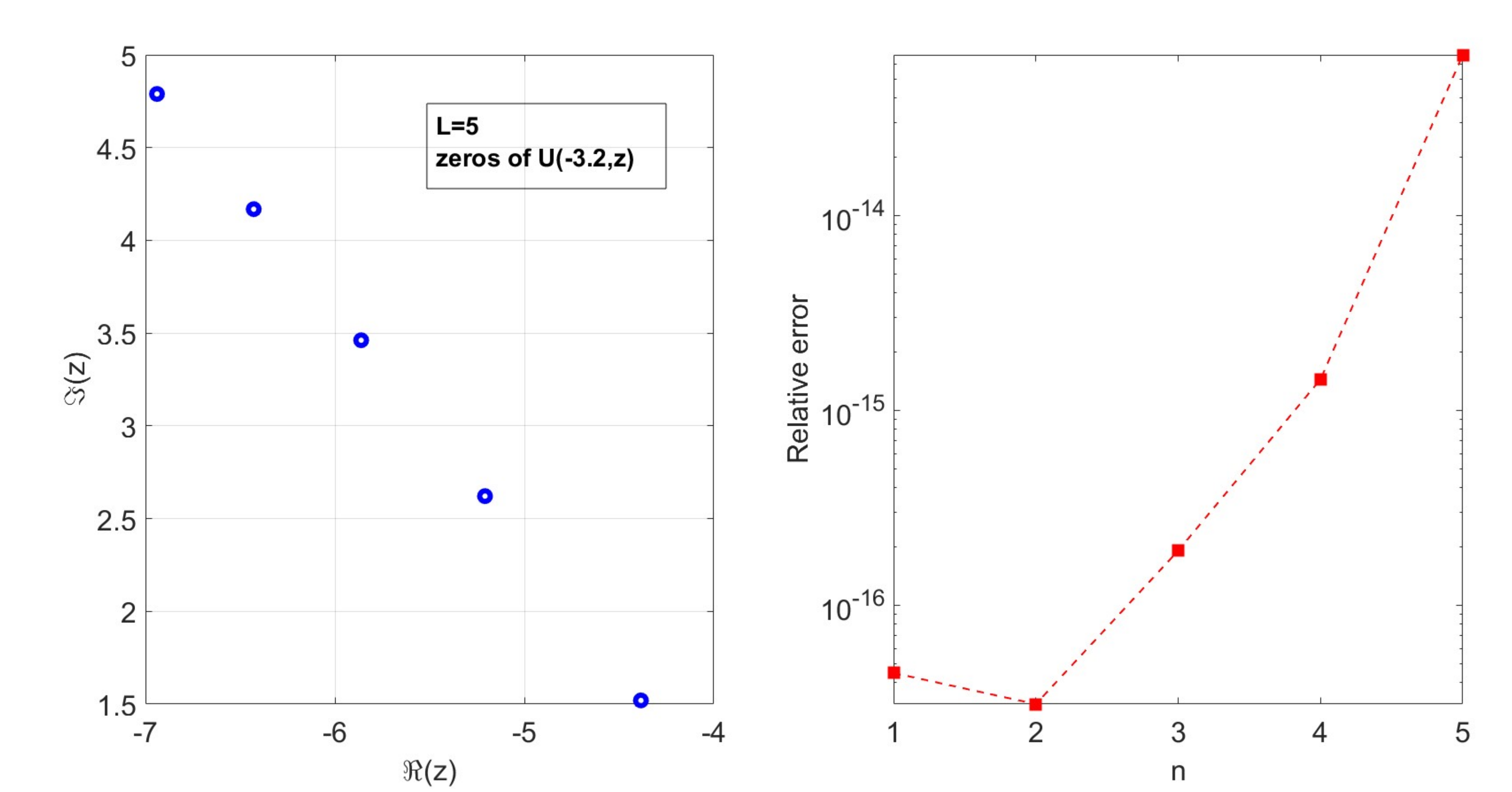}
	\caption{Left: Zeros obtained with the function  \texttt{zerosUaz(a,L)} for $a=-3.2$ and $L=5$. Right: Estimated relative errors obtained. }
	\label{fig:1}
\end{figure}

\begin{figure}[h!]
		\includegraphics[width=1.2\textwidth]{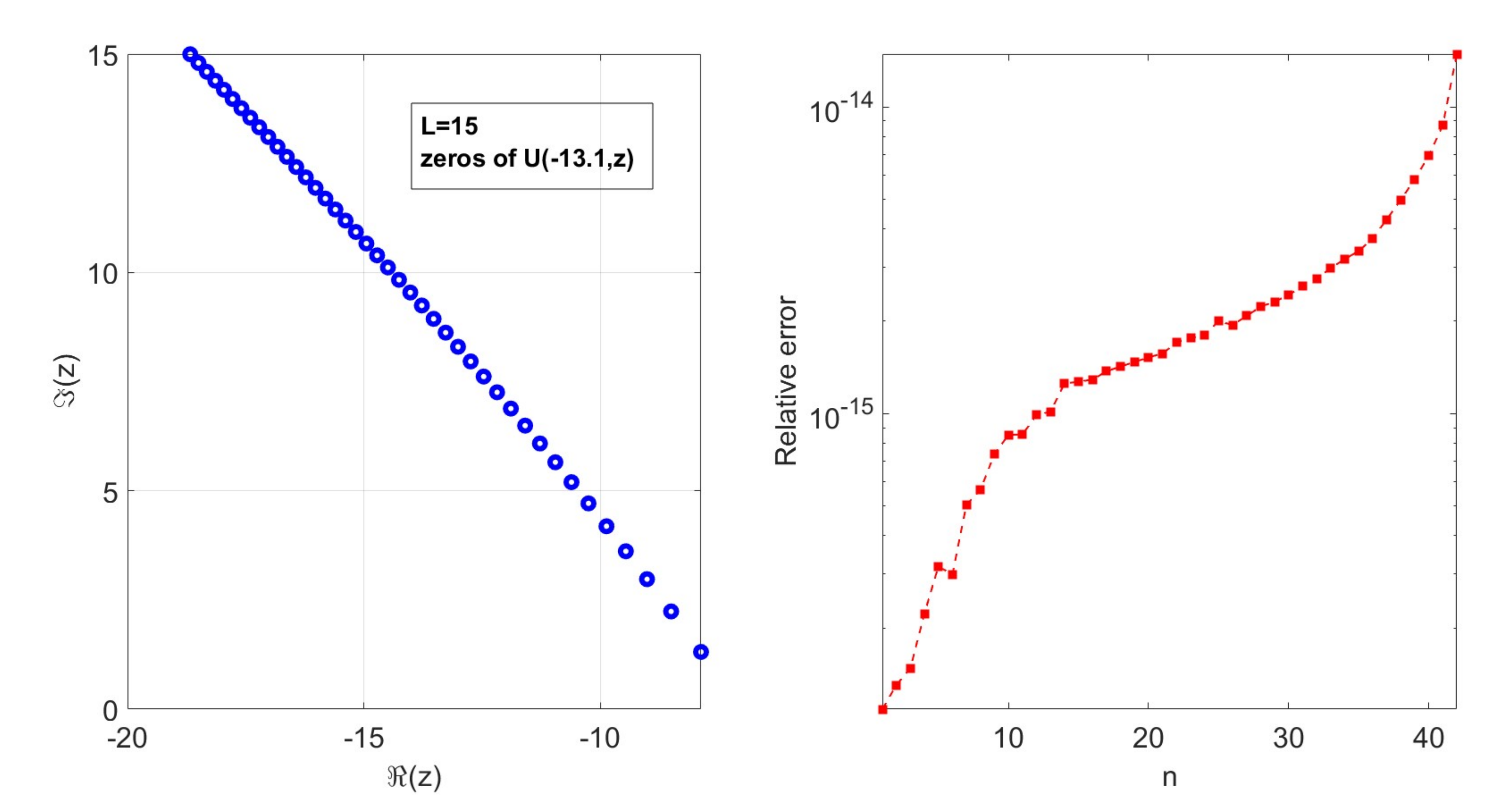}
	\caption{Left: Zeros obtained with the function   \texttt{zerosUaz(a,L)}  for $a=-13.1$ and $L=15$. Right: Estimated relative errors obtained.  }
	\label{fig:2}
\end{figure}

\begin{figure}[h!]
		\includegraphics[width=1.2\textwidth]{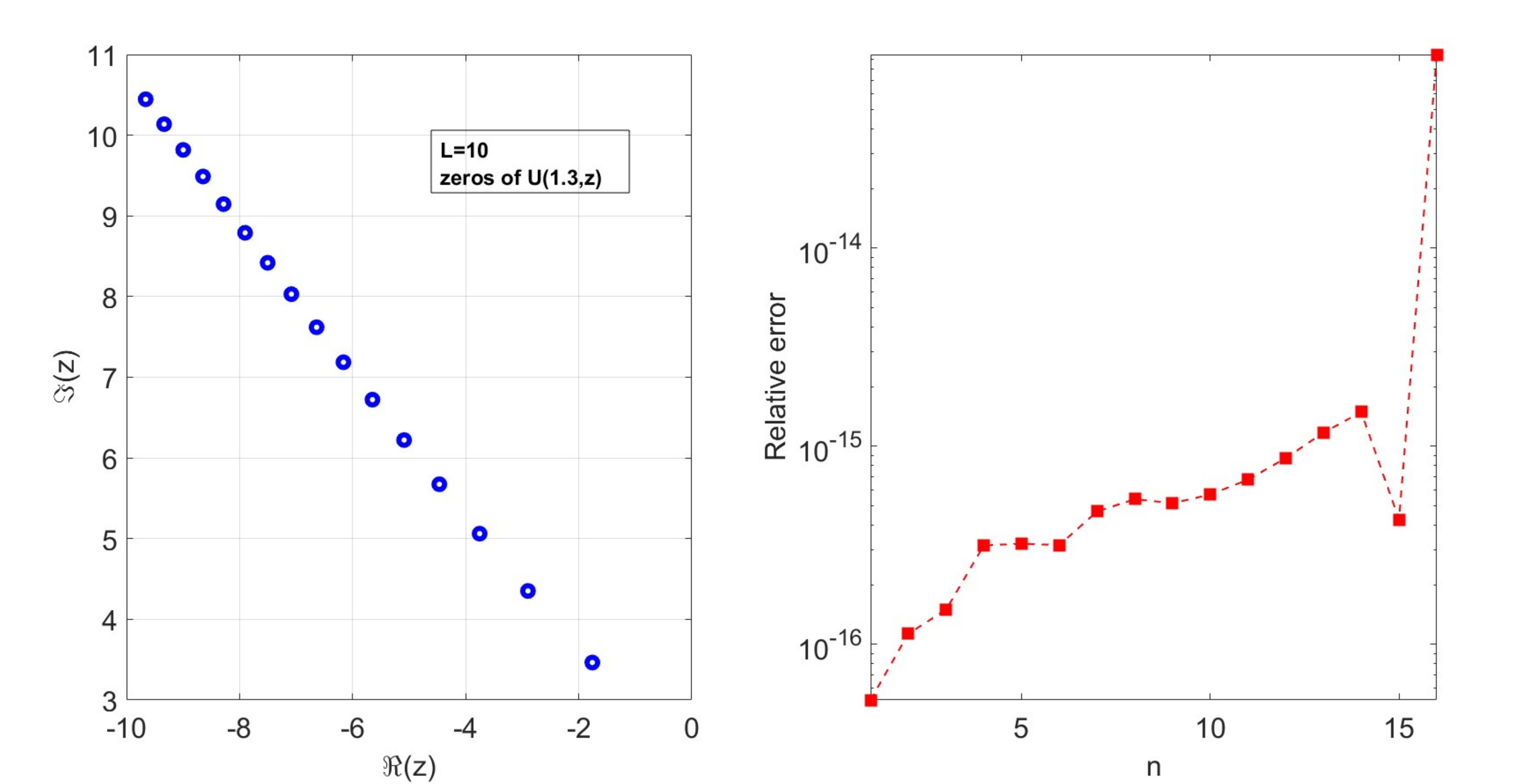}
	\caption{Left: Zeros obtained with   the function   \texttt{zerosUaz(a,L)}  for $a=1.3$ and $L=10$. Right: Estimated relative errors obtained. }
	\label{fig:3}
\end{figure}

\begin{figure}[h!]
		\includegraphics[width=1.2\textwidth]{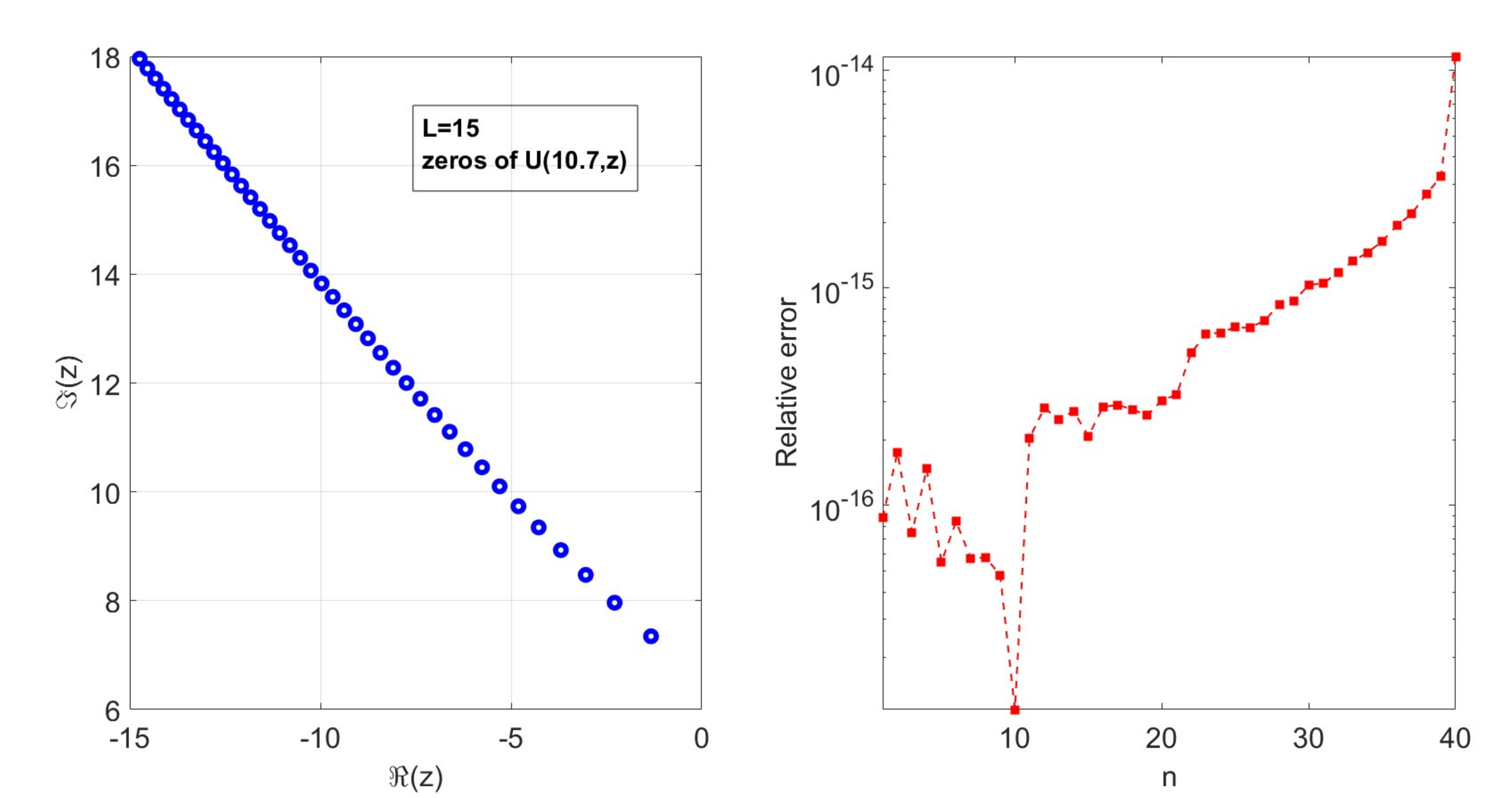}
	\caption{Left: Zeros obtained with  the function   \texttt{zerosUaz(a,L)}  for $a=10.7$ and $L=15$. Right: Estimated relative errors obtained.  }
	\label{fig:4}
\end{figure}

\begin{figure}[h!]
		\includegraphics[width=1.2\textwidth]{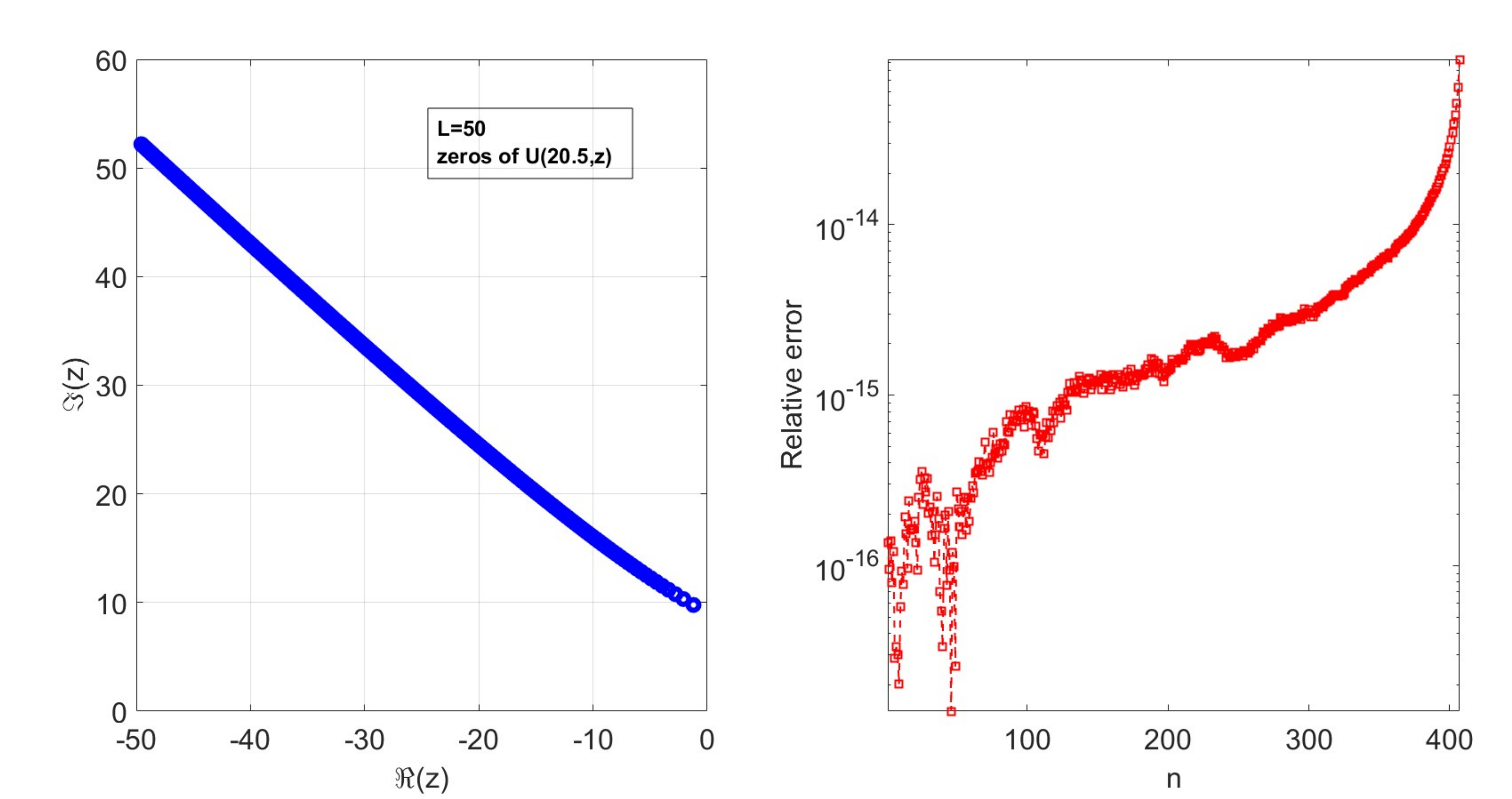}
	\caption{Left: Zeros obtained with  the function   \texttt{zerosUaz(a,L)}  for $a=20.5$ and $L=50$. Right: Estimated relative errors obtained.  }
	\label{fig:5}
\end{figure}

In Figures \ref{fig:1} and \ref{fig:2}, we show two examples of the accuracy obtained with  Algorithm \ref{algo} to calculate the complex zeros 
of $U(a,z)$ ($a<0$)  in the domain of the complex plane $\Im z \in [0,\,L]$, $\Re z<0$. 
Figures \ref{fig:3}, \ref{fig:4} and  \ref{fig:5}  show  examples obtained with the algorithm for positive values of the parameter 
$a$ (function \verb|zerosUazpos(a,L)|)  in the domain of the complex plane $\Re z \in [-L,\,0]$, $\Im z>0$. 
 In Figures  \ref{fig:1} and \ref{fig:3}, where small values of the parameter 
$a$ are considered, the asymptotic approximations
for the computation of the first zero
 have been complemented with few fixed-point iterations.
As can be seen in the figures, almost all zeros are calculated with accuracy significantly better than $\epsilon=10^{-14}$. In the case of the zero closest to the real axis (the one calculated with the greatest estimated relative error), the accuracy is only slightly above this value. This holds true even when the number of zeros to be calculated is very high, as shown, for example, in Figure \ref{fig:5}. In that case, the number of calculated zeros in the domain is $N_z=407$. As an additional check of the accuracy of the zeros computed with our algorithm, we have compared them with Maple values for the zeros (computed with 50 digits) and found that the relative errors obtained are consistent with the accuracy estimated using the inverse of the condition number. We give two examples:
first, the last three zeros computed with our algorithm in the example shown in Figure \ref{fig:2} were (in the order in which they were obtained)
 $z_a=$\verb|-9.008392235290984e+00 + 2.976766819022779e+00i|,
$z_b=$\verb|-8.498829407276800e+00 + 2.237094348893690e+00i| and \\
 $z_c=$\verb|-7.866459577089755e+00 + 1.309795045190692e+00i|. \\
The corresponding zero values computed  with Maple were\\
 $\hat{z}_a=$\verb|-9.00839223529104890...+2.976766819022788250...*I|, \\
$\hat{z}_b=$\verb| -8.4988294072768770...+2.237094348893700886...*I| and \\
  $\hat{z}_c=$\verb|-7.86645957708986376...+1.309795045190640585...*I|. \\
The comparisons give the following relative errors: 
$\epsilon_a=6.99 \times 10^{-15}$, $\epsilon_b=9.05 \times 10^{-15}$ and  $\epsilon_c=1.5 \times 10^{-14}$, respectively.
On the other hand,  the last three zeros computed with our algorithm in the example shown in Figure \ref{fig:5} were
$z_a=$\verb|  -2.784978156368416e+00 + 1.077249953770989e+01i |, \\
$z_b=$\verb|  -2.082060971609168e+00 + 1.031935184018148e+01i  | and \\
 $z_c=$\verb|  -1.204905397657948e+00 + 9.772189846956170e+00i |. \\
The corresponding zero values computed  with Maple were\\
 $\hat{z}_a=$\verb|-2.78497815636791519...+10.77249953770959430...*I  |, \\
$\hat{z}_b=$\verb|-2.08206097160857784... +10.31935184018114877...*I| and \\
$\hat{z}_c=$\verb|-1.20490539765712604... +9.772189846955761085...*I|. \\
The comparisons give the following relative errors: 
$\epsilon_a= 5.23 \times 10^{-14}$, $\epsilon_b=6.41 \times 10^{-14}$ and  $\epsilon_c=9.31 \times 10^{-14}$, respectively.

Regarding the computational efficiency of the algorithm, some typical computation times are shown in Table \ref{table1}
for different values of the parameter $a$ and $L$.
The calculations have been performed in Matlab R2024b, in a computer under Windows 11 (64-bit), processor Intel(R) Core(TM) i5-10210U @ 1.60GHz   2.11 GHz. 
 As can be seen, even a very high number of zeros $N_z$ can be calculated with a low computational cost, demonstrating the efficiency of the numerical scheme.
It is important to mention that the use of Taylor series to compute most fixed-point iterations is largely responsible for the excellent computational efficiency of the algorithm. Although it would have been possible to calculate the quotients \eqref{quot} using the algorithm for computing the function and its derivative, this approach is slower.

In conclusion, the numerical tests conducted demonstrate the accuracy and efficiency of the numerical scheme (implemented in Matlab) we have developed to compute
 zeros of parabolic cylinder functions in domains of the complex plane. As mentioned in the Introduction, these strategies can be extended to
the computation of zeros of other functions that are solutions of second-order ODEs.

\begin{table}[h!]
$$
\begin{array}{cccc}
\hline
   a &  L & N_z & \mbox{CPU time} (s) \\ 
  \hline
        -1.7 & 12 & 23 & 0.023\\
       -1.7 & 60 & 573 & 0.068 \\
      -1.7 & 180 & 5157 & 0..317 \\
      -30.2 & 12 & 31 & 0.018 \\
      -30.2 & 60 &587 & 0.068 \\
      -30.2 & 180 & 5171 &  0.268 \\
           2.3 &  10 &   16          &     0.018  \\
           2.3 &  50&   398      &  0.048     \\
  2.3 &  140&   3120      &  0.228     \\
           20.5 &  10 &   21     &   0.017     \\
           20.5 &  50 &   407      &  0.042           \\
    20.5 &  140 &   3129      &  0.219           \\
\hline
     \end{array}
$$
\caption{\label{table1} Typical CPU times spent by the algorithm. $N_z$ is the number of zeros
calculated in the domain.
}
\end{table}

\section{Declarations}

\subsection{Competing interests}
The authors have no conflict of interest to declare that are relevant to the content of
 this article.

\subsection{Funding}

Financial support was received from {\emph{Ministerio de Econom\'{\i}a y
Competitividad}}, project PID2021-127252NB-I00 (MCIN/AEI/10.13039/501100011033/ FEDER, UE.

\subsection{Acknowledgements}
The authors thank the referees for their helpful suggestions.

\bibliographystyle{spmpsci}
\bibliography{biblio}

\end{document}